\documentclass[twoside, 10pt]{article}

\usepackage{amssymb,amsfonts, amsmath, amsthm}
\usepackage[all,arc]{xy}
\usepackage{graphicx, enumerate}
\usepackage{mathrsfs, datetime}
\usepackage{mathtools, comment, verbatim}
\usepackage{tabularx,tabulary}
\usepackage[top=2cm, bottom=3cm, left=3cm, right=3cm]{geometry}
\usepackage{authblk}
\title{Two-Player Pebbling on Diameter 2 Graphs}
\author{Garth Isaak}
\affil{Department of Mathematics, Lehigh University}
\author{Matthew Prudente}
\affil{Department of Mathematics, Saint Vincent College}
\date{January 19, 2018}                     
\newtheorem{thm}{Theorem}[section]
\newtheorem{lem}[thm]{Lemma}
\newtheorem{prop}[thm]{Proposition}
\newtheorem{cor}[thm]{Corollary}
\newtheorem{fact}[thm]{Fact}
\newtheorem{rslt}[thm]{Result}

\theoremstyle{definition}
\newtheorem{dfn}[thm]{Definition}
\newtheorem{remark}{Remark}

\newtheorem{conj}[thm]{Conjecture}

\newtheorem{shade}{ }

\DeclarePairedDelimiter{\floor}{\lfloor}{\rfloor}

\newcommand{\bthm}{\begin{thm}}
\newcommand{\ethm}{\end{thm}}
\newcommand{\blem}{\begin{lem}}
\newcommand{\elem}{\end{lem}}
\newcommand{\bcor}{\begin{cor}}
\newcommand{\ecor}{\end{cor}}
\newcommand{\bdfn}{\begin{dfn}}
\newcommand{\edfn}{\end{dfn}}
\newcommand{\bconj}{\begin{conj}}
\newcommand{\econj}{\end{conj}}
\newcommand{\bprop}{\begin{prop}}
\newcommand{\eprop}{\end{prop}}
\newcommand{\brslt}{\begin{rslt}}
\newcommand{\erslt}{\end{rslt}}
\newcommand{\bprof}{\begin{proof}}
\newcommand{\eprof}{\end{proof}}
\newcommand{\ds}{\displaystyle}
\newcommand{\bfig}{\begin{figure}}
\newcommand{\efig}{\end{figure}}
\newcommand{\bfact}{\begin{fact}}
\newcommand{\efact}{\end{fact}}
\newcommand{\bshd}{\begin{shade}}
\newcommand{\eshd}{\end{shade}}
\newcommand{\brem}{\begin{remark}}
\newcommand{\erem}{\end{remark}}

\newcommand{\up}{^{\prime}}

\newcommand{\Cprime}{C^{\prime}}

\newcommand{\be}{\begin{enumerate}}
\newcommand{\ee}{\end{enumerate}}
\newcommand{\bi}{\begin{itemize}}
\newcommand{\ei}{\end{itemize}}
\newcommand{\bfl}{\begin{flushleft}}
\newcommand{\efl}{\end{flushleft}}
\newcommand{\bc}{\begin{center}}
\newcommand{\ec}{\end{center}}
\newcommand{\bflr}{\begin{flushright}}
\newcommand{\eflr}{\end{flushright}}
\newcommand{\bas}{\begin{align*}}
\newcommand{\eas}{\end{align*}}

\begin{document}

\maketitle

\begin{abstract}

A \emph{pebbling move} refers to the act of removing two pebbles from one vertex and placing one pebble on an adjacent vertex.  The goal of graph pebbling is: Given an initial distribution of pebbles, use pebbling moves to reach a specified goal vertex called the \emph{root}.  The \emph{pebbling number} of a graph $\pi(G)$ is the minimum number of pebbles needed so every distribution of $\pi(G)$ pebbles can reach every choice of the root.  We introduce a new variant of graph pebbling, a game between two players. One player aims to move a pebble to the root and the other player aims to prevent this.  We  show configurations of various classes of graphs for which each player has a winning strategy.  We will characterize the winning player for a specific class of diameter two graphs.  
\end{abstract}

\section{Introduction}
Graph pebbling can be thought of as an optimization problem where a utility such as gas, electricity, or computing power travels across a network.  While traveling through the network, some amount of the utility may be lost.  A natural question that arises is what is the minimum amount of the utility that is needed to travel the network and arrive at a destination.  

Graph pebbling was originally developed to solve a number theory conjecture posed by Erd\"os \cite{Lemke}.  The goal is to use pebbling moves to place one pebble on a specified vertex $r$ called the \emph{root}.  A \emph{pebbling move} refers to the act of removing two pebbles from one vertex and placing one pebble on an adjacent vertex.  The basic question is, given an initial arrangement of the pebbles called a configuration $C$, can we set one pebble on $r$ through a sequence of pebbling moves.  If so, then $C$ is \emph{r-solvable} \cite{Surglenn}.  We say $C(v)$ is the number of pebbles at vertex $v$ and the \emph{size} of a configuration $C$ is ${\ds \sum_{v \in G} C(v) }$.  We can see in Figure \ref{pm} that if the root is the leftmost vertex, then $C$ is $r$-solvable.  Define $\pi(G,r)$ as the minimum number $m$ such that every configuration of $m$ pebbles is $r$-solvable for a given root $r$. The \emph{pebbling number} ${\ds \pi(G) = \max_{r \in V(G)} \pi(G,r)}$ is the minimum number $t$ such that every configuration of size $t$ is $r$-solvable for every choice of $r$ in $G$.  Graph pebbling is well studied and has numerous variations \cite{hurlgeneral, Surglenn}.

\bfig[!ht]
  \centering
    \includegraphics[scale=1]{triangleline4.mps}
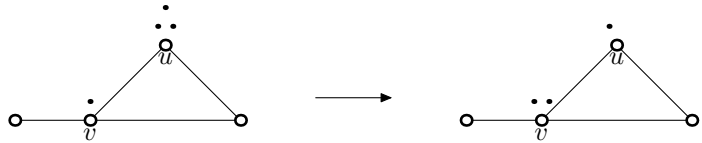
  \caption{An example of a pebbling move on $G$ from $u$ to $v$.}
  \label{pm}
\efig

From this point, all graphs will be connected, finite and simple (no loops or multiedges).  Let $V(G)$ be the set of vertices of $G$ and $|V(G)|$ be the number of vertices in $G$.  The \emph{diameter} of a graph, $diam(G)$, is the longest of all shortest paths in $G$.  The \emph{open neighborhood} of $v$, $N(v)$, is the set of vertices adjacent to but not including $v$.  Likewise, the \emph{closed neighborhood} of $v$, $N[v]$, is the set of vertices adjacent to and including $v$.  Similarly we define the neighborhoods $N(T)$ and $N[T]$ for a set $T$ of vertices as the union of the vertex neighborhoods. Given $S \subseteq V(G)$ and $v \in V(G)$, we say the \emph{S-restricted neighborhood of v}, $N_S(v) = N(v) \cap S$, is the set of neighbors of $v$ contained only in $S$.  Given a graph $G$, the complement $G^{\prime}$ is the graph such that $V(G) = V(G^{\prime})$ and $uv \in E(G^{\prime}) \iff uv \notin E(G)$. We use $K_n$ to denote the complete graph on $n$ vertices.

We introduce a new variation that extends pebbling to a two-person game.

\bdfn
\label{tppdef}
 Let $G$ be a connected graph with specified root $r$ and $C$ a configuration on the vertices of $G$. There are two players, \textit{Mover} and \textit{Defender}. 
 We say a \emph{round} consists of two pebbling moves; the initial move made by Mover and second move made by Defender.  A \emph{turn} will be an individual player's pebbling move.  The \emph{Two-Player Pebbling Game} is then as follows.

\begin{enumerate}
\item Play proceeds in rounds, with Mover pebbling first, then Defender. Each player must take their turn.
\item If Mover pebbles from $u$ to $v$, then Defender can not pebble from $v$ to $u$ in the same round.
\item If $C^{\prime}(r) > 0$ at any time, then Mover wins.
\item If $C^{\prime}(r) = 0$ and there are no more pebbling moves, then Defender wins.
\end{enumerate}

\edfn

Without Rule 2, Mover has limited options to win as otherwise Defender could `undo' any move when a second pebble is placed on a vertex. Without this rule, whenever the graph induced by 
vertices at distance at least 2 from the root is nontrivial there will be arbitrarily large configurations for which Defender wins. If all pebbles are on vertices of nontrivial components then Defender can always pebble to a vertex in the component and never place a second pebble on a vertex in $N(r)$. If Mover places a second pebble on a vertex in $N(r)$ then Defender can pebble back along the same edge. So Mover can never place a pebble on the root. 

Still considering the possibility of a version with no Rule 2, the only interesting cases are when the graph induced by $V - N[r]$ is trivial. A subset of these graphs is what we consider with Rule 2 in Section \ref{sec2}. For the version with no Rule 2 it is not difficult to show that we would get the same conclusions as Theorem \ref{moverwinconfig} except that Defender always has a winning strategy in the exceptional case that $k$ is even, $C_T = k+2$ and exactly one vertex in $T$ is even.

However, if we consider a variant where Rule 1 was relaxed to allow Defender to forfeit their turn and we ignore Rule 2, then Defender always has a winning strategy.

We examine conditions which ensure a win for Mover or a win for Defender.  To do this, we must study how each player will play the game.

\bdfn
A \emph{strategy} for either player is a choice function $\mathcal{S}: \mathcal{C} \to \mathcal{P}$ from the set of all possible configurations $\mathcal{C}$ to the list of all allowed pebbling moves $\mathcal{P}$.
\edfn

By this, of course, we mean a strategy is a method of playing the game based on the possible outcomes of any move.  
Moves for Defender must take into account Mover's previous moves by Rule 2 in Definition \ref{tppdef}.

\bdfn
A strategy $\mathcal{S}$ is \emph{winning} for Mover (or Defender) on a configuration $C$ provided Mover (or Defender) wins playing according to $\mathcal{S}$. We sometimes say `Mover wins' (or `Defender wins') to refer to Mover (Defender) having a winning strategy. 
\edfn

Now we can introduce the values for two-player pebbling.

\bdfn
For a graph $G$ with root $r$, the \emph{rooted-two-player pebbling number}, $\eta(G,r)$, is the minimum number $m$ such that given any configuration of $m$ pebbles, Mover has a winning strategy.  From this, we say the \emph{two-player pebbling number} is ${\ds \eta(G) = \max_{r \in V} \eta(G,r)}$, the minimum number $t$ such that for every configuration of size $t$ and every choice of $r$, Mover has a winning strategy.  However, if for a graph $G$, a root $r$, and arbitrarily large $m$, there exists a configuration of size at least $m\up$, for $m\up > m$, for which Defender has a winning strategy, then we say that $\eta(G,r) = \infty$ and if $\eta(G,r) = \infty$ for some $r$ then $\eta(G) = \infty$.
\edfn

\section{Preliminary Results}

We begin with some basic statements about $\eta(G)$.

\bprop
\label{pieta}
$|V(G)| \leq \pi(G) \leq \eta(G)$.
\eprop

\bprof
$|V(G)| \leq \pi(G)$ is well known as placing a single pebble on all vertices except the root is a distribution that cannot reach the root. For the two-player pebbling number, the moves of a winning strategy suffice for classical pebbling.
\eprof

Notice that if Defender is not forced to pebble in a winning pebbling move sequence for classical pebbling, then equality fails.

The proof for Proposition \ref{n-1proof} is essentially the same as that for the corresponding result for classical pebbling, \cite{Chung}.

\bprop
\label{n-1proof}
If $deg(r) = |V(G)| - 1$, then $\eta(G, r) = |V(G)|$.
\eprop

\bprof
Let $r$ be a vertex with degree $|V(G)| - 1$.  Suppose we have $|V(G)|-1$ pebbles.  If every non-root vertex has 1 pebble, then Defender wins.  So suppose we have $|V(G)|$ pebbles. If we have a configuration with 1 pebble on $r$, then Mover wins.  Suppose we have a configuration with no pebbles on the root.  Then there must exist at least one vertex with at least 2 pebbles on it.  Since Mover begins the game, they will pebble to the root.
\eprof

From this, we get a corollary about the complete graph on $n$ vertices, $K_n$.

\bcor
\label{completeproof}
$\eta(K_n) = n$.
\ecor

\subsection{Sufficient Condition for Infinite $\eta$} 

In this section, we give a simple condition for Defender to have a winning strategy for configurations with an arbitrarily large number of pebbles.  Informally this suggest that `most' graphs will have $\eta(G) = \infty$. However we will see interesting results with finite $\eta(G)$ on certain  structured classes of graphs. 

\bthm
\label{defwin}
For a graph $G$, let $S$ be a cut set of $G$ and let $G_0, G_1, \dots G_k$ be the components of $G-S$ with $r \in G_0$. If for every $v \in S, \, |N(v) - (V(G_0) \cup S)| \geq 2$ and for every $x \in N(S) - (V(G_0) \cup S), \,\, |N(x)-S| \geq 2$, then $\eta(G,r) = \infty$.
\ethm

\bprof
Let $G$ be described as above.  Let $m$ be an arbitrary natural number and $\mathcal{C}$ be the family of configurations with $m$ pebbles on the vertices of $G_1 \cup \cdots \cup G_k$ and no pebbles on $V(G_0) \cup S$. We will show that Defender has a strategy which prevents either player from pebbling from $S$ to $G_0$. Hence no pebble will be placed on the root. With a finite number of pebbles Defender will win.

 Suppose Mover puts a second pebble on a vertex $v \in S$.  Because $|N(v) - (V(G_0) \cup S)| \geq 2$, Defender can pebble to another vertex in $N(v) - (V(G_0) \cup S)$. Thus Defender will never be forced to pebble to $G_0$. If Defender avoids placing a second pebble on a vertex in $S$ then Mover can never pebble to $G_0$.  Let $x \in N(S) - (V(G_0) \cup S)$ and suppose Defender must pebble from $x$.  Because $|N(x)-S| \geq 2$, Defender can pebble from $x$ to a vertex in $N(x)-S$.  Therefore, Defender is never forced to place a second pebble on a vertex in $S$.
\eprof

\bfig[!ht]
  \centering
    \includegraphics[scale=1]{NecForDef1.mps}
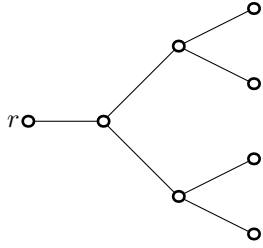
    \caption{A small example for Theorem \ref{defwin}.}
    \label{NecForDef}
\efig

Note that Figure \ref{NecForDef} satisfies the conditions for Theorem \ref{defwin}.  We see that Figure \ref{NecForDef} is a tree, and thus bipartite.  Therefore, trees and bipartite graphs can have an infinite two-player pebbling number. The graph in Figure \ref{smalldiam} has diameter 2.  Thus, a graph $G$ having diameter 2 is not a sufficient condition for a finite value of $\eta(G)$, whereas diameter-2 graphs have classical pebbling number of at most $|V(G)| + 1$ \cite{onpebblinggraphs}.

\bfig[!ht]
  \centering
    \includegraphics[scale=1]{NecForDef2.mps}
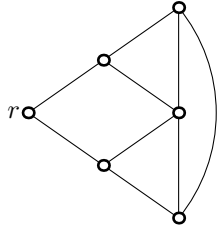
    \caption{A graph with diameter 2 for Theorem \ref{defwin}.}
    \label{smalldiam}
\efig

Observe that grids, $P_n \square P_m$ for $m, n \geq 4$ have infinite $\eta$ because they satisfy the conditions for Theorem \ref{defwin}.  Consider Figure \ref{cartprod}.  

\bfig[!ht]
  \centering
    \includegraphics[scale=1]{pathproduct1.mps}
    \caption{$P_4 \square P_4$.}
  \label{cartprod}
 \efig

It is easy to verify that $\eta(P_4)$ is finite, and hence $\eta(P_4) \cdot \eta(P_4)$ is finite, but $\eta(P_4 \square P_4) = \infty$.  This is in direct contrast to Graham's Conjecture \cite{Chung}, a well studied problem in classical pebbling  which states $\pi(G \square H) \leq \pi(G)\cdot \pi(H)$ for any choice of $G$ and $H$.   So even for a simple Cartesian product of graphs, a two-player pebbling analog of Graham's Conjecture will not hold.  

With Theorem \ref{defwin} we can also quickly see that two-player-pebbling numbers are not monotone with respect to edge additions and deletions. Let $G_1 = K_6$, $G_2$ be the graph in Figure \ref{smalldiam} and $G_3 = P_6$, the path on 6 vertices. We have $E(G_1) \supset E(G_2) \supset E(G_3)$ with $\eta(G_1) = 6$ by Corollary \ref{completeproof},
$\eta(G_2) = \infty$ by Theorem \ref{defwin} and $\eta(G_3) = 35$ as can readily be checked \cite{Me}. 
 However, removing an edge incident to the root can only help Defender. It is straightforward to check that 
 if $e$ is an edge incident to the root then $\eta(G-e,r) \geq \eta(G,r)$ as any strategy that wins for defender will not have Defender pebbling to the root so the same strategy suffices for $G-e$.

\section{Certain Diameter 2 Graphs} \label{sec2}

We move on to the study of two-player pebbling on certain graphs of diameter 2.  Specifically, we characterize the winning player for nearly every configuration for a specific class of diameter 2 graph, characterize the winning player for every configuration as well as $\eta(G)$ for complete multipartite graphs with all part sizes at least 3.

For any two graphs $H$ and $G$ on disjoint vertex sets, the \emph{join} of $H$ and $G$, $H \vee G$, is the graph that contains all edges in $H$, all edges in $G$, and edges connecting every vertex in $H$ with every vertex with $G$.
We use $H \cup G$ to denote the disjoint union of graphs $H$ and $G$. 

Now, we define a subset of diameter 2 graphs.

\bdfn
For $s \geq 1, t \geq 1$ let  $\mathcal{G}_{s,t}$ denote the class of all graphs of the form $ \big( (K_1 \cup K_t^{\prime}) \vee H\big)$  where $H$ is any graph on $s$ vertices.. The root will be the vertex of the $K_1$. We call $S = V(H)$ and $T = V(K_t^\prime)$ with $|S| = s$ and $|T| = t$. 
\edfn

We will examine $\mathcal{G}_{s,t}$ for $t \geq 2$. The special case  $\mathcal{G}_{s,1}$ is more complicated and will be considered in another paper. 

 Figure \ref{kmulti} gives us an example of a graph in $\mathcal{G}_{s,t}$.

\bfig[!ht]
  \centering
    \includegraphics[scale=0.8]{rst1.mps}
    \caption{The class $\mathcal{G}_{s,t}$.}
  \label{kmulti}
\efig

If a starting configuration has two pebbles on any vertex in $S$, then Mover can pebble to the root and win.  So, we say a \emph{non-trivial configuration} on the vertices of $G$ will have 0 or 1 pebbles on vertices in $S$. Let $k$ be the number of vertices in $S$ that are pebble-free.  We say a vertex $v$ is \emph{even} or \emph{odd} corresponding to the parity of
$C(v)$.  A vertex is \emph{pebbled} provided it has at least 1 pebble on it and is \emph{unpebbled} or \emph{pebble-free} otherwise.

We develop a condition on the distribution of pebbles on $T$ based on the pebble-free vertices in $S$.  Informally, it appears that we can compare how many pebbling moves are in $T$ to the number of pebble-free vertices in $S$.  If there are many more pebbling moves in $T$ than pebble-free vertices in $S$, then Mover wins.  Both players must pebble to $S$.  Eventually, $S$ will have no pebble-free vertices and it will be Defender's turn.  They will pebble to $S$; Mover will pebble to $r$ on their next turn.  On the other hand, if there are many more pebble-free vertices than pebbling moves in $T$, Defender wins.  Defender will always have a pebble-free vertex in $S$ to pebble to. We would like a way to count the number of pebbling moves in $T$.  Notice for any vertex $v \in T$ that ${\ds \floor[\Big]{\frac{C(v)}{2}} }$ will tell us the number of pebbling moves on $v$.  We have the following definition.

\bdfn
We say ${\ds C_T = \sum_{v \in T} \floor[\Big]{\frac{C(v)}{2}}}$ is the number pebbling moves in $T$ with configuration $C$.  

We will use $k$ to denote the number of pebble-free vertices in $S$ for a given configuration. 
\edfn

We will see that if there are $k$ pebble-free vertices in $S$ and $C_T \geq k + 3$, then Mover has a winning strategy.  If $C_T \leq k$, then Defender has a winning strategy.  If $C_T = k + 2$ or $k + 1$, then it depends on the parity of $k$ and the structure of $S$ to find the winning player.

\subsection{When $k$ is odd} 

Lemma \ref{moverconfiguration} is the base case for induction when $k$ is odd.

\blem
\label{moverconfiguration}
Let $G \in \mathcal{G}_{s,t}$, $t \geq 2$ and $C$ be a non-trivial configuration with 1 pebble-free vertex in $S$.  Mover has a winning strategy if and only if ${\ds C_T \geq 2}$.
\elem

\bprof
Suppose $C_T \geq 2$.  Mover will pebble to the unpebbled vertex.  Now there is one more move in $T$ and all vertices in $S$ have a pebble on them.  Defender must pebble to a vertex in $S$, placing a second pebble on a vertex.  Mover pebbles to $r$ and wins.

Conversely, suppose $C_T \leq 1$.  If $C_T = 0$, then there are no pebbling moves in $T$ and Defender wins.  Suppose $C_T = 1$.  Since there is 1 pebbling move in $T$, all the vertices in $T$ without the pebbling move have 0 or 1 pebble on them.  Mover has two choices, to pebble to the unpebbled vertex in $S$ or to place a second pebble on a vertex in $S$.  If Mover pebbles to the pebble-free vertex, then for the new configuration $\Cprime$, $\Cprime_T = 0$.  There are no more pebbling moves and Defender wins.  So suppose Mover pebbles to a pebbled vertex in $S$.  If they can, then Defender will pebble to a pebble-free vertex in $S$ or $T$ and win. If all vertices in $T$ are pebbled, then Defender will place a second pebble on one vertex $v$ in $T$, yielding an extra pebbling move. Obverse that Mover must now pebble from $v$ leaving it pebble-free.  Mover and Defender have the same options as earlier. Thus either Defender wins or pebbles to a pebbled vertex $w \not = v$ in $T$. Defender cannot pebble to $v$ as they cannot pebble back on the same edge just used by Mover. Again there is an extra pebbling move, now from $w$. Again, Mover and Defender have the same options as earlier. However, now $v$ is pebble free and Defender can pebble to $v$ and win.  
\eprof

\blem
\label{oddfloormover}
Let $G \in \mathcal{G}_{s,t}$, $t \geq 2$ and $C$ be a non-trivial configuration with $k$ pebble-free vertices in $S$.  If $k$ is odd and ${\ds C_T \geq k + 1}$, then Mover has a winning strategy on $G$.
\elem

\bprof By induction on $k$. \newline

\emph{Base:} Lemma \ref{moverconfiguration}.

 \emph{Induction:}
Let $k$ be odd and $C_T \geq k+1$.  Mover will pebble to a pebble-free vertex in $S$.  If Defender places a second pebble on a vertex in $S$, Mover wins.  If Defender pebbles to a pebble-free vertex in $S$, then there are $k-2$ pebble-free vertices in $S$ and the resulting configuration $\Cprime$ has $\Cprime_T = C_T - 2$.  Thus $\Cprime_T \geq  k-1 = (k-2)+1$.  Hence, by induction, Mover has a winning strategy.
\eprof

Next is a result when Defender has a winning strategy. This does not depend on the parity of $k$. Later the bound will be improved for $k$ even.

\blem
\label{oddfloordef}
Let $G \in \mathcal{G}_{s,t}$, $t \geq 2$ and $C$ be a non-trivial configuration with $k$ pebble-free vertices in $S$.  If  ${\ds C_T \leq k}$, then Defender has a winning strategy on $G$.
\elem

\bprof By induction on the number of pebbles. \newline 

\emph{Base:} If $C_T = 0$ there are no pebbling moves and Defender wins. This includes the cases with 0 or 1 pebbles.

\emph{Induction:} If $k = 0$ then $C_T = 0$ and the base applies.	If $k = 1$,  then either $C_T = 0$, and the base applies, or $C_T = 1$ and by Lemma \ref{moverconfiguration}, Defender wins.	So assume $k \geq 2$.		If Mover moves to a pebble-free vertex in $S$ then, as $k \geq 2$ a pebble-free vertex will remain and Defender will also pebble to a pebble-free vertex in $S$.	The resulting configuration $\Cprime$ has fewer pebbles, $k-2$ pebble-free vertices in $S$ and $\Cprime_T = C_T - 2 \leq k - 2$.  If Mover places a second pebble on a vertex in $S$ then Defender will pebble from that vertex back to a vertex in $T$.  The resulting configuration $\Cprime$ has fewer pebbles, $k+1$ pebble-free vertices and  $\Cprime_T \leq  C_T  \leq k < k+1$.  In each case, by induction, Defender has a winning strategy. \eprof

So for $k$ odd, we have the following:

\begin{table}[!h]
\centering
\begin{tabular}{|c|c|} 
\hline
 Initital Value of $C_T$  & Winning Player   \tabularnewline \hline \hline 
 $C_T \geq k + 1$    &     Mover       \tabularnewline \hline
 $C_T \leq k$   &     Defender        \tabularnewline\hline 
\end{tabular}
\caption{Value of $C_T$ and its winning player for $k$ odd}
\label{koddwinningplayer}
\end{table}

\subsection{When $k$ is even}

When the number of pebble-free vertices on $S$ is even things are more difficult as one case will depend on the details of the structure of the graph on $S$. We first get bounds that hold in general.  

\blem
\label{evenfloormover}
Let $G \in \mathcal{G}_{s,t}$, $t \geq 2$ and $C$ be a non-trivial configuration with $k$ pebble-free vertices in $S$.  If $k$ is even and ${\ds C_T \geq k + 3}$, then Mover has a winning strategy.
\elem

\bprof
By induction on $k$. \newline

\emph{Base:} Let $k = 0$ and $C_T \geq 3$.  Mover will pebble to $S$, placing a second pebble on one of the vertices.  Defender will pebble back to $T$ or pebble to the root or to a pebbled vertex in $S$ and lose.  The new configuration $\Cprime$ has $\Cprime_T \geq 2$ and now $k = 1$.  By Lemma \ref{moverconfiguration}, Mover wins.

\emph{Induction:} Let $C_T \geq k+ 3$ for $k \geq 1$. Mover will pebble to a pebble-free vertex.  If Defender places a second pebble on a vertex in $S$, then Mover wins.  If Defender pebbles to a pebble-free vertex in $S$, then the new configuration $\Cprime$ has $\Cprime_T = C_T - 2 \geq k + 3 - 2 = k + 1 = (k-2)+3$.  Since $S$ now has $k-2$ pebble-free vertices, Mover has a wining strategy by induction.
\eprof

We will forgo the case when $C_T = k + 2$ for now and leave it for its own section. 

For the next lemma, observe that we have already covered the case  ${\ds C_T \leq k}$ in Lemma \ref{oddfloordef}.

\blem
\label{evenfloordefender1}
Let $G \in \mathcal{G}_{s,t}$, $t \geq 2$ and $C$ be a non-trivial configuration with $k$ pebble-free vertices in $S$.  If $k$ is even and ${\ds C_T \leq k + 1}$, then Defender has a winning strategy.
\elem

\bprof
By induction on $k$. \newline

\emph{Base:} Let $k = 0$.  If $C_T = 0$, then $T$ has no pebbling move and Defender wins.	If $C_T = 1$, then all but one vertex in $T$ has at most 1 pebble on it.	If Mover moves to a pebble-free vertex then no moves remain and Defender wins.  If Mover places a second pebble on a vertex in $S$ then Defender will pebble from that vertex back to a vertex in $T$.  The resulting configuration $\Cprime$ has $k=1$ pebble-free vertex and  $\Cprime_T \leq  C_T  \leq 1$.
Defender has a winning strategy by Lemma \ref{moverconfiguration}.

\emph{Induction:} Let $k$ be even and $C_T \leq k + 1$.  If Mover pebbles to a pebble-free vertex in $S$, then Defender will as well.  The new configuration $\Cprime$ has $k-2$ pebble-free vertices and $\Cprime_T = C_T - 2 \leq k -1 = (k-2)+1$.  By induction, Defender has a winning strategy.  If Mover places a second pebble on a vertex in $S$, Defender will pebble from that vertex to a vertex in $T$.  The resulting configuration $C^{\prime\prime}$ has $k + 1$ pebble-free vertices in $S$ and  $C^{\prime\prime}_T \leq C_T \leq k + 1$. Defender has a winning strategy by Lemma \ref{oddfloordef}.
\eprof

So for $k$ even, we have the following:

\begin{table}[!h]
\centering
\begin{tabular}{|c|c|} 
\hline
 Initital Value of $C_T$  & Winning Player   \tabularnewline \hline \hline 
 $C_T \geq k + 3$    &     Mover       \tabularnewline \hline
 $C_T \leq k + 1$   &     Defender        \tabularnewline\hline 
\end{tabular}
\caption{Value of $C_T$ and its winning player for $k$ even}
\label{kevenwinningplayer}
\end{table}

\subsection{When $k$ is even and $C_T = k+2$} 

The case $k$ even and $C_T = k + 2$ is more difficult to evaluate. The configuration of $S$ and which vertices in $S$ are pebbled as well as how many vertices in $T$ have a non-zero even number of pebbles on them will determine which player has a winning strategy.  Each player's strategy changes a little.  Mover's goal is to force Defender to pebble to a vertex in $T$ with an odd number of pebbles on it.  This will increase the number of pebbling moves in $T$ and yield one of Mover's winning configurations described in an earlier section.  Defender will try to pebble to a vertex in $T$ with an even number of pebbles on it.  This adds no new pebbling moves and yields one of Defender's winning configurations from an earlier section.

In this section we consider some cases determined by the parity of the number of pebbles on vertices in $T$. The next section will consider the structure of $S$ and the number of pebbles on a particular vertex in $T$. 

\blem
\label{evenfloormover2}
Let $G \in \mathcal{G}_{s,t}$, $t \geq 2$ and $C$ be a non-trivial configuration with $k$ pebble-free vertices in $S$.  If $k$ is even and ${\ds C_T = k + 2}$ and for all $v \in T$, $C(v)$ is odd, then Mover has a winning strategy.
\elem

\bprof
By induction on $k$. \newline

\emph{Base:} Let $k = 0$ and $C_T = 2$ with every vertex in $T$ having an odd number of pebbles on it.  Mover will pebble to $S$, placing a second pebble on one of the vertices.  Defender will pebble back to $T$ or lose.  Since every vertex in $T$ has an odd number of pebbles, the new configuration $\Cprime$ has $\Cprime_T = 2$ with 1 unpebbled vertex in $S$.  By Lemma \ref{moverconfiguration} Mover wins.

\emph{Induction:} Let $k$ be even and $C_T = k + 2$ for $k \geq 1$. Mover will pebble to a free vertex.  If Defender places a second pebble on a vertex in $S$, then Mover wins.  If Defender pebbles to a pebble-free vertex in $S$, the new configuration $\Cprime$ has $\Cprime_T = C_T - 2 = (k + 2) - 2 = (k-2) + 2$.  Since $S$ now has $k-2$ pebble-free vertices, Mover has a wining strategy by induction.
\eprof

Now, we look at the case when some vertices in $T$ have an even number of pebbles on them.  This becomes more difficult.  The strategies for each player depends on how many pebbles on are the vertex with an even number of pebbles.

\blem
\label{evenfloordefx0}
Let $G \in \mathcal{G}_{s,t}$, $t \geq 2$ and $C$ be a non-trivial configuration with $k$ pebble-free vertices in $S$.  If $k$ is even and ${\ds C_T  = k + 2}$ and there is either at least one $x \in T$  such that $C(x) = 0$ or at least two vertices $x, y \in T$  such that $C(x)$ and $C(y)$ are even, then Defender has a winning strategy.
\elem

\bprof

By induction on $k$. \newline

\emph{Base:}  Let $k = 0$ and $C_T = 2$.  Mover must place a second pebble on a vertex in $S$.  Defender will pebble from that vertex to the pebble-free vertex in $T$ or to an even vertex in $T$.  For the new configuration $\Cprime$, we have $\Cprime_T = 1$ and $k = 1$.  Thus by Lemma \ref{moverconfiguration}, Defender wins.

\emph{Induction:} Let $k$ be even and $C_T \geq k + 2$.  Mover can place a second pebble on a vertex in $S$ or pebble to a pebble-free vertex in $S$.  If Mover places a second pebble on a vertex in $S$, then Defender will pebble to the unpebbled vertex in $T$ or to an even vertex in $T$, not adding any pebbling moves to $T$.  For our new configuration $\Cprime$, we have $\Cprime_T = k + 1$ and $k$ is now odd.  Hence, Defender wins by Lemma \ref{oddfloordef}.  If Mover pebbles to a pebble-free vertex in $S$, then Defender will also pebble to a pebble-free vertex in $S$.  Now for our new configuration $\Cprime$, we have $\Cprime_T = k$ and there are $k-2$ pebble-free vertices in $S$.  Since there were no pebbling moves back to $T$, we can see that $T$ will still have at least one pebble-free vertex or at least two even vertices.  Thus, Defender wins by induction.
\eprof

So for $k$ even and $C_T = k+2$, we have the following:

\begin{table}[!h]
\centering
\begin{tabular}{|c|c|c|} 
\hline
Number of Even Vertices in $T$  & Initial Value of $C_T$ & Winning Player   \tabularnewline \hline \hline 
None    &  $C_T = k+2$ &    Mover       \tabularnewline \hline
At least one pebble-free or two even   &   $C_T = k+2$  &   Defender        \tabularnewline\hline 
\end{tabular}
\caption{Number of even vertices in $T$ and its winning player when $k$ is even and $C_T = k+2$.}
\label{kevennumbereven}
\end{table}

\subsection{A New Game} 

One case remains,  $k$ even, $C_T = k + 2$ and exactly one vertex in $T$ has a positive, even number of pebbles. For this case we will introduce a new game, the element selecting game that will be equivalent to determining the winner of the pebbling game.
The element selecting game will allow us to determine the winner for a few more cases as well as when $S$ is a complete multipartite graph. However, other cases will remain were we cannot give a simple characterization for which player will win and the game illustrates why these cases are particularly difficult.

We will refer to this configuration frequently enough that it is convenient to give it a name.

\bdfn
Let $G\in \mathcal{G}_{s,t}$, $t \geq 2$. If $C$ is a non-trivial configuration with $k$ pebble-free vertices in $S$ with the following properties we call $C$ a \emph{boundary configuration}. 
The configuration has $k$ even and ${\ds C_T  = k + 2}$.  There is one even vertex $x \in T$  such that $C(x) \geq 2$ and all other vertices in $T$ have an odd number of pebbles.
\edfn

To motivate the element selecting game we first give two more pebbling lemmas.  

\blem
\label{oneeven2mover}
Let $G\in \mathcal{G}_{s,t}$, $t \geq 2$ with a boundary configuration $C$ having $k$ pebble-free vertices. 
Mover has a winning strategy if there exists a pebbled vertex $v \in S$ such that all its neighbors in $S$ are pebbled. 
\elem

\bprof
Mover will pebble from $x$ to $v$.  Defender can either pebble to a neighbor of $v$ or pebble to an odd vertex in $T$. If Defender pebbles to a neighbor of $v$, then that vertex will have two pebbles on it and Mover wins. If Defender pebbles to an odd vertex in $T$, then they will add a pebbling move.  Now our new configuration $\Cprime$ has $k + 1$ pebble-free vertices in $S$ and $\Cprime_T = k + 2$.  By Lemma \ref{oddfloormover}, Mover has a winning strategy.
\eprof

  If $S$ is independent, then the conditions for Lemma \ref{oneeven2mover} will hold vacuously.  Here is a configuration for Defender's winning strategy.

\blem
\label{oneeven2defender}
Let $G\in \mathcal{G}_{s,t}$, $t \geq 2$ with a boundary configuration $C$ having $k$ pebble-free vertices. 
Defender has a winning strategy if $C(x) = 2$ and for every pebbled vertex $v \in S$, there exists at least one pebble-free neighbor in $u \in S$.  
\elem

\bprof
Mover can pebble to a pebbled vertex or an unpebbled vertex.  If Mover pebbles to a pebbled vertex $v$, then Defender will pebble from $v$ to its pebble-free neighbor, which exists by our hypothesis.  Now $k$ is unchanged and our new configuration $\Cprime$ is such that $\Cprime_T = k + 1$.  By Lemma \ref{evenfloordefender1}, Defender has a winning strategy.  If Mover pebbles to an unpebbled vertex, then Defender will pebble from $x$ to another vertex in $S$ which is pebble-free, which exists because $k$ is even and at least 2.  By Lemma \ref{evenfloordefx0}, Defender has a winning strategy.
\eprof

The goal of Mover is to have a pebbled closed neighborhood for some vertex in $S$ and still have at least 2 pebbles on the one even vertex in $T$ and the goal of Defender is to prevent this. Thus we will look at sets corresponding to neighborhoods
of vertices in $S$. 

\bdfn \label{esg1}
Let $N_1, N_2, \dots N_k$ be a collection of  subsets, possibly empty and/or intersecting, from a universal set $U$.  There are two players, Mary and Dan.  Each player will take turns, Mary beginning and Dan following, selecting one element from $U$.  The play continues for up to $j$ rounds.  Mary wins if at least one of the subsets $N_i$ has every one of its elements selected and Dan wins if none of the $N_i$'s has been completely selected by the end of the $j^{th}$ round. If there exists a subset $N_m$ which is empty, then we say Mary wins vacuously.  We call this the \emph{Element Selecting Game with $j$ rounds}. 
\edfn 

This game directly relates to two-player pebbling with the boundary configurations.

\bdfn 
Let $S_0$ be the pebble-free vertices of $S$ and $S_1$ be the pebbled vertices of $S$. 
\edfn

\bdfn
\label{newgame}
Let $G\in \mathcal{G}_{s,t}$, $t \geq 2$ with a boundary configuration $C$ having $k = 2j$ pebble-free vertices. 
Define $\mathcal{E}(G,C,j)$ as the instance of the Element Selecting Game with $j$ rounds constructed in the following way: Let $U=S_0$, the set of unpebbled vertices in $S$.  For every vertex $v_i \in S$, let $N_i = N[v_i] \cap U$.  
\edfn

Now, we can show that the two games are equivalent when we restrict Two-Player Pebbling to this current case.

\bthm
\label{twogamesequal}
Let $G\in \mathcal{G}_{s,t}$, $t \geq 2$ with a boundary configuration $C$ having $k = 2j$ pebble-free vertices. 
 Let $\mathcal{E}(G,C,j)$ be the instance of the Element Selecting Game with $j$ rounds constructed from $G$ as in Definition \ref{newgame}. Mover has a winning strategy in the Two-Player Pebbling Game if and only if Mary has a winning strategy for the Element Selecting Game.
\ethm

\bprof 
We can assume that Defender will never place a 2nd pebble on a vertex in $S$ as in that case Mover wins. 
Mover will always pebble from a vertex with an odd number of pebbles in $T$ and Defender will always pebble from $x$. 
The choices of vertices pebbled to in $S$ and elements selected will be made according to the strategy in the corresponding game. Mary wins the Element Selecting Game if and only if the corresponding Two-Player Pebbling Game reaches a state where
Lemma \ref{oneeven2mover} applies and Dan wins if and only if the corresponding Two-Player Pebbling Game reaches a state where Lemma \ref{oneeven2defender} applies. 
\eprof

Theorem \ref{twogamesequal} and the Element Selecting Game provide short alternate proofs of Lemmas \ref{oneeven2mover} and \ref{oneeven2defender}.
The following Corollary is also straightforward from the sets in the Element Selecting Game $\mathcal{E}(G,C,j)$.

\bcor
\label{oneeven2moversuff}
Let $G\in \mathcal{G}_{s,t}$, $t \geq 2$ with a boundary configuration $C$ having $k$ pebble-free vertices. 
Mover has a winning strategy if $C(x) \geq k + 2$.
\ecor

\bprof
Every set in $\mathcal{E}(G,C,j)$ has size at most $k$. 
\eprof

So for $k$ even with $C_T = k+2$ and one even vertex $x \in T$, we have the following:

\begin{table}[h!]
\centering
\begin{tabular}{|p{8cm}|c|c|c|} 
\hline
Structure of $S$ 					 &  $C(x)$ 	& Initial value of $C_T$	 & Winning Player   \tabularnewline \hline \hline 
Any structure    						 &  $C(x) \geq k+2$  &  $C_T = k+2$ 	 &  Mover       \tabularnewline \hline
Some pebbled vertex with all pebbled neighbors    	 & $ C(x) \geq 2 $	 &  $C_T = k+2$ 	 &  Mover       \tabularnewline \hline
All pebbled vertices have an unpebbled neighbor    & $ C(x) = 2 $	 &  $C_T = k+2$ 	 &  Defender       \tabularnewline \hline

\end{tabular}
\caption{Structure of $S$ and its Winning Player when $k$ is even, $C_T = k+2$ and there is one even vertex $x \in T$. Cases not covered are equivalent to the Element Selecting Game.}
\label{kevenwinningplayer}
\end{table}

We can continue to use the Element Selecting Game with 2 rounds to characterize configurations for the case $C(x) = 4$. The proof is omitted but straightforward. 

\bcor
\label{oneeven4}
Let $G\in \mathcal{G}_{s,t}$, $t \geq 2$ with a boundary configuration $C$ having $k$ pebble-free vertices with $C(x) = 4$.
Mover has a winning strategy if  there exists a vertex $v$ in $S_0$ that for every vertex $u \in S_0$ that either:
\be[a)]
\item there is some vertex $w \in S_1$ such that $N_{S_0}(w) = \{v\}$ or $\{u, v\}$, or
\item $N_{S_0}(u) = \{v\}$.  
\ee

Defender has a winning strategy if for every vertex $v$ in $S_0$ there exists a vertex $u \in S_0$ such that there is no vertex $w \in S_1$ with $N_{S_0}(w) = \{v\}$ or $\{u, v\}$ and (b) $N_{S_0}(u) \neq \{v\}$.
\ecor

The condition for Mover to win extends to $C(x) \geq 4$ while the condition for Defender to win requires $C(x) \leq 4$. It seems unlikely that we will get easily stated complete results when 
$C(x) \geq 6$. Instead we can say something for the case when $S$ is highly structured.

While we will be able to determine $\eta(G)$ for complete multipartite graphs with part sizes at least 3 using the information we already have, in order to have a complete characterization of the winner for all configurations on this class we can use the Element Selecting Game. 

\bthm
Let $G\in \mathcal{G}_{s,t}$ with a boundary configuration $C$ having $k$ pebble-free vertices with $C(x) \leq k$. Assume 
the graph induced by $S$ is a complete multiparitite graph with $m$ partite sets having
$0 \leq k_1 \leq \cdots \leq k_m$ unpebbled vertices respectively. When $k_m \geq \frac{k}{2}$ Mover has a winning strategy
if and only if $C(x) \geq 2(k - k_m) + 2$. When $k_m < \frac{k}{2}$ Mover has a winning strategy if and only if 
$C(x) \geq k+2$.
\ethm

\bprof Define $A_i$ as the unpebbled vertices from the $i^{th}$ part of the graph induced by $S$.  In each case Mover has a winning strategy if $C(x) \geq k+2$ by Corollary \ref{oneeven2moversuff}. So we may assume that $C(x) \leq k+1$ and by Theorem \ref{twogamesequal} consider the  corresponding Element Selecting Game from 
Definition \ref{newgame} with sets ${\ds N_i = \bigcup_{i \not = j} A_i}$ or ${\ds N_i = \{v \} \bigcup_{i \not = j} A_i }$ where the vertices in the parts are $A_i$ depending on whether or not $v \in A_j$ is pebbled.

If $k_m < \frac{k}{2}$ Mover has a winning strategy if  
$C(x) \geq k+2$ by Corollary \ref{oneeven2moversuff}. When $C(x) \leq k+1$ the corresponding Element Selecting Game from 
Definition \ref{newgame} has strictly less than $\frac{k}{2}$ rounds and sets ${\ds N_i = \bigcup_{i \not = j} A_i }$ or ${\ds N_i = \{v \} \bigcup_{i \not = j} A_i}$ where the vertices in the parts ar $A_i$ depending on whether or not $v \in A_j$ is pebbled. If Defender always selects from $A_i$ with the most unselected elements the new values $k'_i$ for unselected elements (i.e., un pebbled vertices) maintain the property that the largest is at most the sum of the others. With less than $\frac{k}{2}$ rounds at least 2 elements remain and by the property there are at least two distinct $A_i$ with unselected elements
and hence no $N_i$ has been selected and Defender wins. 

If $k_m > \frac{k}{2}$ and $C(x) \geq 2(k - k_m) + 2$ there are $(k - k_m)$ rounds and as the smallest $N_i$ has size
at most $k - k_m$ Mover has a winning strategy by picking elements from the smallest $N_i$. If $C(x) < 2(k - k_m) + 2$ there are $j < (k - k_m)$ rounds and at the end there are at least $2k_m - k + 2$ unselected elements. If Defender always selects 
from $A_i$ with the most unselected elements, they can guarantee that $A_m$ has at most $k_m - (k - k_m - 1) = 2k_m + 1$ unselected elements and hence if the new $k_m'$ is largest it has an unselected element and some other $A_i$ has an unselected element. If Mover's selections force some other $k'_i$ to be largest then as at least two elements remain at least two $A_i$ are nonempty as in the previous paragraph. In each case Defender wins. 
\eprof

\subsection{Determining $\eta(\mathcal{G}_{s,t}, r)$}

Now we summarize results from previous lemmas and proceed to determine $\eta(G,r)$ for $G \in \mathcal{G}_{s,t}$.

\bthm
\label{moverwinconfig}
Let $G$ in $\mathcal{G}_{s,t}$, $t \geq 2$ and $C$ be a configuration with $k$ pebble-free vertices in $S$.  Then we have the following: \newline

\begin{center} 
\begin{tabular}{|p{6cm}|p{6cm}|} 
\hline
Mover has a winning strategy on $G$ provided				 &  Defender has a winning strategy on $G$ provided   \tabularnewline \hline \hline 
$k$ is odd and $C_T \geq k+1 $   						 &   $k$ is odd and $C_T \leq k$ 	      \tabularnewline \hline
$k$ is even and $C_T \geq k+3$  					  	 & $k$ is even and $C_T \leq k+1$  	       \tabularnewline \hline
$k$ is even and $C_T = k+2$ and all vertices in $T$ are odd          & $k$ is even and $C_T = k+2$ and $T$ has atleast one unpebbled vertex or two even vertices	     \tabularnewline \hline
\end{tabular}
\end{center}

And if $k$ is even and $C_T = k+2$ and exactly one vertex in $T$ is even, then the game is equivalent to the Element Selecting Game.

\ethm

Obtaining the winning configurations for Mover allow us to get $\eta(G ,r)$ for $G \in \mathcal{G}_{s,t}$.

\bthm
\label{GinG}
If $G \in \mathcal{G}_{s,t}$, $t \geq 2$,  then
$\eta(G,r) = 
\begin{cases}
t + 2s + 4, & s \text{ is even}\\
t + 2s + 3, & s \text{ is odd.}
\end{cases}$
\ethm

\bprof

\smallskip
In each case we refer to Theorem \ref{moverwinconfig} to determine the winner. 

\emph{Case 1:} Let $s$ be even.  A configuration of $t + 2s + 3$ pebbles on the vertices of $G$ which gives Defender a winning strategy is the following: in $T$, leave one vertex pebble-free, put one pebble on $t-2$ vertices and the remaining $2s + 5$ pebbles on one vertex and keep $S$ pebble-free.  With this configuration, $C_T = s + 2$ with one vertex in $T$ having no pebbles on it and Defender wins.

Now suppose there are $m \geq t + 2s + 4$ pebbles on the vertices in $G$.  Let $k$ of the vertices in $S$ be pebble-free.  Thus there are $(s-k)$ pebbles in $S$. Note that $k \leq s$. Now there are $m - (s-k) \geq t + 2s + 4 - s + k = t+ s + k + 4$ pebbles on the vertices in $T$. If all vertices in $T$ are odd then  $C_T = \left\lfloor \frac{s + k + 4}{2}  \right\rfloor \geq k+2$ and Mover wins. 
If at least one vertex in $T$ is even then $C_T > \left\lfloor \frac{s + k + 4}{2} \right\rfloor \geq k+2$ and Mover wins.

\smallskip

\emph{Case 2:} Let $s$ be odd.  A configuration of $t + 2s + 2$ pebbles on the vertices of $G$ which gives Defender a winning strategy is the following: place 1 pebble on any vertex in $S$, place 1 pebble on $t-1$ vertices in $T$ and the remaining $2s + 1$ pebbles on the other vertex in $T$. With this configuration, $C_T = s$ and there are $s-1$ pebble-free vertices in $S$, with $s-1$ even and  Defender wins.

Now suppose there are $m \geq t + 2s + 3$ pebbles on the vertices in $G$.  Let $k$ of the vertices in $S$ be pebble-free.  Thus there are $(s-k)$ pebbles in $S$.  Now there are $m - (s-k) \geq t + 2s + 3 - s + k = t+ s + k + 3$ pebbles on the vertices in $T$. If $k$ is odd then $C_T = \left\lfloor \frac{s + k + 3}{2} \right\rfloor \geq k+1$ and Mover wins.
If $k$ is even then $s \geq k+1$ since they have different parity. 
If all vertices in $T$ are odd then $C_T = \left\lfloor \frac{s + k + 3}{2} \right\rfloor \geq \left\lfloor \frac{(k+1) + k + 3}{2} \right\rfloor = k+2$ and Mover wins. 
If at least one vertex in $T$ is even then $C_T > \left\lfloor \frac{s + k + 3}{2} \right\rfloor \geq \left\lfloor \frac{(k+1) + k + 3}{2} \right\rfloor  = k+2$ and Mover wins.
\eprof

Since for complete multipartite graphs we are in the class $\mathcal{G}_{s,t}$ for any choice of root, 
using Theorem \ref{GinG} we can determine $\eta(G)$ when all part sizes are at least 3.

\bthm
\label{CBbiggerthan3}
If $3 \leq a_1 \leq a_2 \leq \dots \leq a_m < n$ and ${\ds \sum_{i=1}^m a_i = n}$, then
$$
\eta(K_{a_1, a_2, \dots, a_m}) = 
\begin{cases}
2n - a_1 + 3, & n-a_1 \text{ is even}\\
2n - a_1 + 2, & n - a_1 \text{ is odd.}
\end{cases}
$$
\ethm

\bprof
If $r$ is in a part of size $a_k$ then we can apply Theorem \ref{GinG} with ${\ds s = n - a_k }$ and $t = a_k - 1$. Then for this choice of $r$,   $\eta(K_{a_1, a_2, \dots, a_m},r)$ is either
$(a_k -1) + 2(n - a_k) + 4 $ or $(a_k - 1) + 2(n - a_k) + 3$ depending on the parity of $s$. It is easy to check that $k = 1$ attains the maximum.  
\eprof

Observe that for the special case of complete bipartite graphs $K_{u,v}$ with $3 \leq u \leq v$ and $n = u + v$ we have 
$\eta(K_{u,v}) = 2n - u + 3 = u + 2v + 3$ when $n - u = v$ is even and $\eta(K_{u,v}) = 2n - u + 2 = u + 2v + 2$ when $n - u = v$ is odd. It is straightforward to check that $\eta(K_1,v) = v + 4 = n + 3$ when $v \geq 3$. It is more complicated but can also be checked that $\eta(K_{2,v}) = v + 7 = n + 5$ when $v \geq 3$. Thus the multipartitie formulas do not extend and we can get a quite different value when at least one part size is 1 or 2.

\section{Conclusion}

Two-player pebbling in some other situations are examined in \cite{Me} as well as other forthcoming papers. We will briefly indicate some of what is known and not known here. 

Theorem \ref{GinG} shows that $\eta(G,r)$ depends only on the parity of $s = |S|$ for $G \in \mathcal{G}_{s,t}$ when $t \geq 2$. For  $G \in \mathcal{G}_{s,1}$, $\eta(G,r)$ can depend on the structure of $S$.  This case will, in part relate to the Element Selecting Game but includes even more complicated strategies for the players. Use $n = s + 2$ for the number of vertices in $G$.  We can show that when $s \geq 4$,  $n = s + 2 \leq \eta(G,r) \leq 2s + 2 = 2n - 2$ for $G \in \mathcal{G}_{s,1}$ and that these bounds are best possible. In particular 
$\eta(G,r) = s + 2 = n$ if $S$ induces an independent set so that $G$ is $K_{2,s}$ with root in the part of size 2 and $\eta(G,r) = 2s + 2 = 2n - 2$ if $S$ induces a clique so that $G$ is a complete graph minus an edge with root one of the ends of the omitted edge. 

While paths, $P_n$ seem like a potentially straightforward situation this appears not to be the case. For path powers $P_n^k$,  $2 \leq k \leq n-5$, $\eta(P_n^k) = \infty$. Also $\eta(P_n^{n-1}) = n$ (the complete graph), 
$\eta(P_n^{n-2}) = 2n - 2$ (the complete graph minus an edge) and $\eta(P_n^{n-3})$, $\eta(P_n^{n-4})$ are both finite.  
For simple paths recall the  classical pebbling number $\pi(P_n) = 2^{n-1}$. We can show that 
$\eta(P_n) \leq \frac{3}{2} \cdot 2^{n-1} - n$ for $n \geq 4$. We can show for small values that this upper bound is not the correct value. Playing the `obvious strategy' for each player is not optimal. For somewhat more complicated strategies that we think may be optimal, we can determine the pebbling number under these fixed strategies for each player. However we do not have an exact value of $\eta(P_n)$.

\bibliographystyle{plain}
\bibliography{myreferences}

\end{document}